\documentclass[12pt]{article}
\usepackage{amssymb}
\usepackage[mathscr]{eucal}
\usepackage{amsmath,amssymb,latexsym,theorem,bbm}
\newfont{\msa}{msam10 scaled\magstep1}
\newfont{\ssmsa}{msam9}
\newfont{\smsa}{msam10}
\newfont{\sms}{msbm10}
\newfont{\sseufb}{eufb9}
\newfont{\seufb}{eufb10}
\newfont{\eufb}{eufb10 scaled\magstep1}
\newfont{\eusb}{eusb10 scaled\magstep1}
\newfont{\hcmr}{cmr17 scaled\magstep5}

\newcommand{\AstkoneKn}{\raise-7mm\hbox{\hcmr*}%
          ^{\hspace*{-5.3mm}\raise2.8mm\hbox{$\scriptstyle K_n$}}%
           _{\hspace*{-6.1mm}\raise2.5mm\hbox{$\scriptstyle k=1$}}}
\newcommand{\Astkonen}{\raise-7mm\hbox{\hcmr*}%
          ^{\hspace*{-5.3mm}\raise2.8mm\hbox{$\scriptstyle n$}}%
           _{\hspace*{-6.1mm}\raise2.5mm\hbox{$\scriptstyle k=1$}}}

\newcommand{\DS}{\displaystyle}

\newcommand{\SC}{\scriptstyle}

\newcommand{\CC}{\mathbb{C}}

\newcommand{\EE}{\mathrm{E}}
\newcommand{\LL}{\mathbb{L}}

\newcommand{\NN}{\mathbb{N}}
\newcommand{\PP}{\mathrm{P}}

\newcommand{\RR}{\mathbb{R}}
\newcommand{\TT}{\mathbb{T}}

\newcommand{\ZZ}{\mathbb{Z}}

\newcommand{\sSS}{\raise-0.5truemm\hbox{\sms S}}

\newcommand{\cI}{\mathcal{I}}

\newcommand{\cP}{\mathcal{P}}

\newcommand{\cN}{\mathcal{N}}

\newcommand{\ww}{\mathrm{w}}
\newcommand{\dd}{\mathrm{d}}
\newcommand{\ee}{\mathrm{e}}
\newcommand{\qq}{\mathrm{q}}
\renewcommand{\Re}{\mathrm{Re}\,}

\newcommand{\hmu}{\widehat{\mu}}
\newcommand{\hdelta}{\widehat{\delta}}
\newcommand{\homega}{\widehat{\omega}}
\newcommand{\hgamma}{\widehat{\gamma}}
\newcommand{\hpi}{\widehat{\pi}}

\newcommand{\hDelta}{\widehat{\Delta}}
\newcommand{\hS}{\widehat{S}}
\newcommand{\hTT}{\widehat{\TT}}
\newcommand{\hRR}{\widehat{\RR}}
\newcommand{\hG}{\widehat{G}}

\newcommand{\sleq}{\mbox{\ssmsa\hspace*{0.1mm}\symbol{54}\hspace*{0.1mm}}}

\renewcommand{\leq}{\mbox{\msa\hspace*{0.9mm}\symbol{54}\hspace*{0.9mm}}}
\renewcommand{\geq}{\mbox{\msa\hspace*{0.9mm}\symbol{62}\hspace*{0.9mm}}}

\newcommand{\bone}{\mathbbm{1}}
\newcommand{\vare}{\varepsilon}

\newcommand{\proofend}{\hfill\mbox{$\Box$}}

\newcommand{\csD}{{\SC\mathcal{D}}}

\newcommand{\distre}{\stackrel{\csD}{=}}

\numberwithin{equation}{section}

\theoremstyle{change} \theorembodyfont{\em}
\newtheorem{Lem}{Lemma.}[section]
\newtheorem{Thm}[Lem]{Theorem.}

\theorembodyfont{\rm}
\newtheorem{Def}[Lem]{Definition.}

\begin{document}

\begin{center}
 {\bfseries\Large Weakly infinitely divisible measures}\\[3mm]
 {\bfseries\Large on some locally compact Abelian groups}\\[5mm]
 {\bfseries\large M\'aty\'as $\text{Barczy}^{*,\diamond}$} {\large and}
 {\bfseries\large Gyula $\text{Pap}^*$}
\end{center}

\vskip0.5cm

* Faculty of Informatics, University of
  Debrecen, Pf. 12, H--4010 Debrecen, Hungary;\\
  E--mail: barczy@inf.unideb.hu (M. Barczy); papgy@inf.unideb.hu (G. Pap).

$\diamond$ Corresponding author.

\vspace*{-1mm}

{\small\textbf{Abstract.} On the torus group, on the group of \
$p$--adic integers and on the \ $p$--adic solenoid we give a
construction of an arbitrary weakly infinitely divisible probability
measure using a random element with values in a product of
(possibly infinitely many) subgroups of \ $\RR.$ \
As a special case of our results, we have a new construction of
the Haar measure on the \ $p$--adic solenoid.}

\indent {\small\textbf{Running head.} Weakly infinitely divisible
measures on some LCA groups.}

\section{Introduction}

Weakly infinitely divisible probability measures play a very
important role in limit theorems of probability theory, see for example, the books
of Parthasarathy \cite{PAR}, Heyer \cite{HEY}, the papers of Bingham \cite{BIN},
Yasuda \cite{YAS}, Barczy, Bendikov and Pap \cite{BARBENPAP},
and the Ph.D. theses of Gaiser \cite{GAI}, Tel\"oken \cite{TEL}, Barczy \cite{BAR}.
They naturally arise as possible limits of triangular arrays described as follows.

Let \ $G$ \ be a locally compact Abelian topological group
having a countable basis of its topology. We also suppose that
\ $G$ \ has the \ $T_0$--property, that is, \ $\bigcap_{U\in\cN_{e}}U=\{e\},$ \
where \ $e$ \ denotes the identity element of \ $G$ \ and
\ $\cN_e$ \ is the collection of all Borel neighbourhoods of \ $e.$ \
(By a Borel neighbourhood \ $U$ \ of \ $e$ \ we mean a Borel subset of \
$G$ \ for which there exists an open subset \ $\widetilde U$ \ of \ $G$ \
 such that \ $e\in \widetilde U\subset U.$) \ Let us consider a probability measure
 \ $\mu$ \ on \ $G$ \ and let \ $\{X_{n,k}:n\in\NN,\,k=1,\dots,K_n\}$ \
 be an array of rowwise independent random elements with values in \ $G$ \
 satisfying the infinitesimality condition
 $$
   \lim_{n\to\infty}\max_{1\sleq k\sleq K_n}\PP(X_{n,k}\in G\setminus U)=0,
     \qquad \forall\; U\in\cN_e.
 $$
If the row sums \ $\sum_{k=1}^{K_n}X_{n,k}$ \ of such an array
converge in distribution to \ $\mu$ \ then \ $\mu$ \ is necessarily weakly infinitely
 divisible, see, e.g., Parthasarathy \cite[Chapter IV, Theorem 5.2]{PAR}.
Moreover, Parthasarathy \cite[Chapter IV, Corollary 7.1]{PAR} gives
a representation of an arbitrary weakly infinitely divisible measure
on \ $G$ \ in terms of a Haar measure, a Dirac measure, a Gauss
measure and a Poisson measure on \ $G$.

In this paper we consider the torus group, the group
of \ $p$--adic integers and the \ $p$--adic solenoid.
For these groups, we give a construction of an arbitrary weakly
infinitely divisible measure using real random variables.
 For each of these three groups, the construction consists in:
\renewcommand{\labelenumi}{{\rm(\roman{enumi})}}
 \begin{enumerate}
  \item finding a group, say \ $G_0,$ \ which is a product of (possibly infinitely many)
        subgroups of \ $\RR.$ \ (We furnish \ $G_0$ \ with the product topology.
        Note that \ $G_0$ \ is not necessarily locally compact.)
  \item finding a continuous homomorphism \ $\varphi: G_0\to G$ \
        such that for each weakly infinitely divisible measure \ $\mu$ \ on \ $G,$ \
        there is a probability measure \ $\mu_0$ \ on \ $G_0$ \ with the property
        \ $\mu_0(\varphi^{-1}(B))=\mu(B)$ \ for all Borel subsets \ $B$ \ of \ $G.$ \
        (The probability measure \ $\mu_0$ \ on \ $G_0$ \ will be given as the distribution
         of an appropriate random element with values in \ $G_0.$)
 \end{enumerate}
 Since \ $\varphi$ \ is a homomorphism, the building blocks of
 \ $\mu$ \ (Haar measure, Dirac measure, Gauss measure and Poisson
 measure) can be handled separately.

We note that, as a special case of our results, we have a new
construction of the Haar measure on the \ $p$-adic integers and the
\ $p$-adic solenoid. Another kind of description of the Haar measure
on the \ $p$-adic integers can also be found in Hewitt and Ross
\cite[p. 220]{HR}. One can find a construction of the Haar measure
on the \ $p$-adic solenoid in Chistyakov \cite[Section 3]{CHI}. It
is based on Hausdorff measures and rather sophisticated, while our
simpler construction (Theorem \ref{TwidS}) is based on a
probabilistic method and reflects the structure of the \ $p$-adic
solenoid.

\section{Parametrization of weakly infinitely divisible measures}
\label{par}

Let \ $\NN$ \ and \ $\ZZ_+$ \ denote the sets of positive and of
nonnegative integers, respectively.
The expression ``a measure \ $\mu$ \ on \ $G$'' means a measure \
$\mu$ \ on the \ $\sigma$--algebra of Borel subsets of \ $G$.
\ The Dirac measure at a point \ $x\in G$ \ will be denoted by \
$\delta_x$.

\begin{Def} \
A probability measure \ $\mu$ \ on \ $G$ \ is called
 \emph{infinitely divisible} if for all \ $n\in\NN,$ \ there exist a
 probability measure \ $\mu_n$ \ on \ $G$ \ such that \ $\mu=\mu_n^{*n}$. \
 The collection of all infinitely divisible measures on
  \ $G$ \ will be denoted by \ $\mathcal{I}(G)$. \
 A probability measure \ $\mu$ \ on \ $G$ \ is called
 \emph{weakly infinitely divisible} if for all \ $n\in\NN,$ \ there exist a
 probability measure \ $\mu_n$ \ on \ $G$ \ and an element \ $x_n\in G$ \ such
 that \ $\mu=\mu_n^{*n}*\delta_{x_n}$. \
 The collection of all weakly infinitely divisible measures on
  \ $G$ \ will be denoted by \ $\mathcal{I}_{\ww}(G)$. \
\end{Def}

Note that \ $\mathcal{I}(G)\subset\mathcal{I}_{\ww}(G),$ \ but in general
\ $\mathcal{I}(G)\ne\mathcal{I}_{\ww}(G).$ \ Clearly,
\ $\mathcal{I}(G)=\mathcal{I}_{\ww}(G)$ \ if all the Dirac measures on
\ $G$ \ are infinitely divisible. In case of the torus and the $p$--adic solenoid,
\ $\mathcal{I}(G)=\mathcal{I}_{\ww}(G),$ \ see Sections \ref{widT} and \ref{widS};
and in case of the $p$--adic integers, \ $\mathcal{I}(G)\ne\mathcal{I}_{\ww}(G),$ \
see the example in Section \ref{widD}. We also remark that Parthasarathy \cite{PAR}
 and Yasuda \cite{YAS} call weakly infinitely divisible measures on
 \ $G$ \ infinitely divisible measures.

We recall the building blocks of weakly infinitely divisible
measures. The main tool for their description is Fourier
transformation. A function \ $\chi:G\to\TT$ \ is said to be a character of
\ $G$ \ if it is a continuous homomorphism, where \ $\TT$ \ is the topological group of
 complex numbers \ $\{\ee^{ix}:-\pi\leq x<\pi\}$ \ under multiplication (for more details on \ $\TT,$ \
 see Section \ref{widT}). The group of all characters of \ $G$ \ is called
 the character group of \ $G$ \ and is denoted by \ $\hG$. \
 The character group \ $\hG$ \ of \ $G$ \ is also a locally compact Abelian $T_0$--topological group
having a countable basis of its topology (see, e.g., Theorems 23.15 and 24.14 in Hewitt and Ross \cite{HR}).
For every bounded measure \ $\mu$ \ on \ $G$, \ let \ $\hmu:\hG\to\CC$ \ be
 defined by
 $$
  \hmu(\chi):=\int_G\chi\,\dd\mu,
  \qquad \chi\in\hG.
 $$
This function \ $\hmu$ \ is called the \emph{Fourier transform} of \
$\mu$. \ The usual properties of the Fourier transformation can be
found, e.g., in
 Heyer \cite[Theorem 1.3.8, Theorem 1.4.2]{HEY}, in Hewitt and Ross
 \cite[Theorem 23.10]{HR} and in Parthasarathy
 \cite[Chapter IV, Theorem 3.3]{PAR}.

If \ $H$ \ is a compact subgroup of \ $G$ \ then \ $\omega_H$ \ will
denote the
 Haar measure on \ $H$ \ (considered as a measure on \ $G$\,) \ normalized by
 the requirement \ $\omega_H(H)=1$.
\ The normalized Haar measures of compact subgroups of \ $G$ \ are
the only
 idempotents in the semigroup of probability measures on \ $G$ \ (see, e.g.,
 Wendel \cite[Theorem 1]{WEN}).
For all \ $\chi\in\hG$,
 \begin{equation}\label{homega}
  \homega_H(\chi)
  =\begin{cases}
    1 & \text{if \ $\chi(x)=1$ \ for all \ $x\in H$,}\\
    0 & \text{otherwise,}
   \end{cases}
 \end{equation}
 i.e., \ $\homega_H=\bone_{H^\perp},$ \ where
 $$
  H^\perp:=\big\{\chi\in\hG:\text{$\chi(x)=1$ \ for all \ $x\in H$}\big\}
 $$
 is the annihilator of \ $H$.
\ Clearly \ $\omega_H\in\cI_\ww(G)$, \ since \
$\omega_H*\omega_H=\omega_H$.

Obviously \ $\delta_x\in\cI_\ww(G)$ \ for all \ $x\in G$.

A \emph{quadratic form} on \ $\hG$ \ is a nonnegative continuous
function \ $\psi:\hG\to\RR_+$ \ such that
 $$
  \psi(\chi_1\chi_2)+\psi(\chi_1\chi_2^{-1})=2(\psi(\chi_1)+\psi(\chi_2))
  \qquad\text{for all \ $\chi_1,\chi_2\in\hG$.}
 $$
The set of all quadratic forms on \ $\hG$ \ will be denoted by \
$\qq_+(\hG)$. \ For a quadratic form \ $\psi\in\qq_+(\hG)$, \ there
exists a unique probability
 measure \ $\gamma_\psi$ \ on \ $G$ \ determined by
 $$
  \hgamma_\psi(\chi)=\ee^{-\psi(\chi)/2}\qquad\text{for all \ $\chi\in\hG$},
 $$
 which is a symmetric Gauss measure (see, e.g., Theorem 5.2.8 in Heyer
 \cite{HEY}).
Obviously \ $\gamma_\psi\in\cI_\ww(G)$, \ since
 \ $\gamma_\psi=\gamma_{\psi/n}^{*n}$ \ for all \ $n\in\NN$.

For a bounded measure \ $\eta$ \ on \ $G$, \ the
 \emph{compound Poisson measure} \ $\ee(\eta)$ \ is the probability measure
 on \ $G$ \ defined by
 $$
  \ee(\eta)
  :=\ee^{-\eta(G)}
    \left(\delta_e+\eta+\frac{\eta*\eta}{2!}+\frac{\eta*\eta*\eta}{3!}
          +\cdots\right),
 $$
 where \ $e$ \ is the identity element of \ $G.$ \
The Fourier transform of a compound Poisson measure \ $\ee(\eta)$ \
is
 \begin{equation}\label{cP}
  (\ee(\eta))\:\widehat{}\:(\chi)
  =\exp\left\{\int_G(\chi(x)-1)\,\dd\eta(x)\right\},
   \qquad\chi\in\hG.
 \end{equation}
Obviously \ $\ee(\eta)\in\cI_\ww(G)$, \ since
 \ $\ee(\eta)=\big(\ee(\eta/n)\big)^{*n}$ \ for all \ $n\in\NN$.

In order to introduce generalized Poisson measures, we
recall the notions of a local inner product and a L\'evy measure.

\begin{Def}\label{lip} \
A continuous function \ $g:G\times\hG\to\RR$ \ is called a
 \emph{local inner product} for \ $G$ \ if
\renewcommand{\labelenumi}{{\rm(\roman{enumi})}}
 \begin{enumerate}
  \item for every compact subset \ $C$ \ of \ $\hG$, \ there exists
         \ $U\in\cN_e$ \ such that
         $$\chi(x)=\ee^{ig(x,\chi)}\qquad
           \text{for all \ $x\in U$,\quad$\chi\in C$,}$$
  \item for all \ $x\in G$ \ and \ $\chi,\chi_1,\chi_2\in\hG$,
         $$g(x,\chi_1\chi_2)=g(x,\chi_1)+g(x,\chi_2),\qquad
           g(-x,\chi)=-g(x,\chi),$$
  \item for every compact subset \ $C$ \ of \ $\hG$,
         $$\sup\limits_{x\in G}\sup\limits_{\chi\in C}|g(x,\chi)|<\infty,
           \qquad
           \lim\limits_{x\to e}\sup\limits_{\chi\in C}\vert g(x,\chi)\vert=0.$$
 \end{enumerate}
\end{Def}

Parthasarathy \cite[Chapter IV, Lemma 5.3]{PAR} proved the existence
of a local inner product for an arbitrary locally compact Abelian \
$T_0$--topological group having a countable basis of its topology.

\begin{Def} \
A measure \ $\eta$ \ on \ $G$ \ with values in \ $[0,+\infty]$ \
is said to be a \emph{L\'evy measure} if \ $\eta(\{e\})=0$,
 \ $\eta(G\setminus U)<\infty$ \ for all \ $U\in\cN_e$, \ and
 \ $\int_G(1-\Re\chi(x))\,\dd\eta(x)<\infty$ \ for all \ $\chi\in\hG$.
\ The set of all L\'evy measures on \ $G$ \ will be denoted
  by \ $\mathbb{L}(G)$.
\end{Def}

We note that for all \ $\chi\in\hG$ \ there exists \ $U\in\cN_e$ \
such that
 \begin{equation}\label{compare}
  \frac{1}{4}g(x,\chi)^2\leq 1-\Re\chi(x)\leq\frac{1}{2}g(x,\chi)^2,
  \qquad x\in U.
 \end{equation}
 Thus the requirement \ $\int_G(1-\Re\chi(x))\,\dd\eta(x)<\infty$ \
 can be replaced by the requirement that
 \ $\int_G g(x,\chi)^2\,\dd\eta(x)<\infty$ \
 for some (and then necessarily for any) local inner product \ $g$.

For a L\'evy measure \ $\eta\in\LL(G)$ \ and for a local inner
product \ $g$
 \ for \ $G$, \ the \emph{generalized Poisson measure} \ $\pi_{\eta,\,g}$
 \ is the probability measure on \ $G$ \ defined by
 $$
  \hpi_{\eta,\,g}(\chi)
  =\exp\left\{\int_G\big(\chi(x)-1-ig(x,\chi)\big)\,\dd\eta(x)\right\}\qquad
  \text{for all \ $\chi\in\hG$}
 $$
 (see, e.g., Chapter IV, Theorem 7.1 in Parthasarathy \cite{PAR}).
Obviously \ $\pi_{\eta,\,g}\in\cI_\ww(G)$, \ since
 \ $\pi_{\eta,\,g}=\pi_{\eta/n,\,g}^{*n}$ \ for all \ $n\in\NN$.
\ Note that for a bounded measure \ $\eta$ \ on \ $G$ \ with
 \ $\eta(\{e\})=0$ \ we have \ $\eta\in\LL(G)$ \ and
 \ $\ee(\eta)=\pi_{\eta,\,g}*\delta_{m_g(\eta)}$, \ where the element
 \ $m_g(\eta)\in G$, \ called the \emph{local mean} of \ $\eta$ \ with
  respect to the local inner product \ $g$, \ is uniquely defined by
 $$
  \chi(m_g(\eta))=\exp\left\{i\int_G g(x,\chi)\,\dd\eta(x)\right\}\qquad
  \text{for all \ $\chi\in\hG$.}
 $$
(The existence of a unique local mean is guaranteed by Pontryagin's
duality theorem.)

Let \ $\cP(G)$ \ be the set of quadruplets \ $(H,a,\psi,\eta)$, \
where \ $H$ \ is a compact subgroup of \ $G$, \ $a\in G$, \
$\psi\in\qq_+(\hG)$ \ and \ $\eta\in\LL(G)$. \ Parthasarathy
\cite[Chapter IV, Corollary 7.1]{PAR} proved the following
 parametrization for weakly infinitely divisible measures on \ $G$.

\begin{Thm}[Parthasarathy]\label{LCA1} \
Let \ $g$ \ be a fixed local inner product for \ $G$. \ If \
$\mu\in\cI_\ww(G)$ \ then there exists a quadruplet
 \ $(H,a,\psi,\eta)\in\cP(G)$ \ such that
 \begin{equation}\label{LH}
  \mu=\omega_H*\delta_a*\gamma_\psi*\pi_{\eta,\,g}.
 \end{equation}
Conversely, if \ $(H,a,\psi,\eta)\in\cP(G)$ \ then
 \ $\omega_H*\delta_a*\gamma_\psi*\pi_{\eta,\,g}\in\cI_\ww(G)$.
\end{Thm}
In general, this parametrization is not one--to--one (see
Parthasarathy \cite[p.112, Remark 3]{PAR}).

We say that \ $\mu\in\mathcal{I}_{\ww}(G)$ \ has a non-degenerate idempotent factor if
\ $\mu=\mu'*\nu$ \ for some probability measures \ $\mu'$ \ and \ $\nu$ \ such that
\ $\nu$ \ is idempotent and \ $\nu\ne\delta_e.$ \
Yasuda \cite[Proposition 1]{YAS} proved the following characterization of weakly
infinitely divisible measures on \ $G$ \ without non-degenerate idempotent factors, i.e.,
weakly infinitely divisible measures on \ $G$ \ for which in the representation \eqref{LH} the
compact subgroup \ $H$ \ is \ $\{e\}.$ \

\begin{Thm}[Yasuda]\label{THM_YASUDA}
  A probability measure \ $\mu$ \ on \ $G$ \ is weakly infinitely divisible
  without non-degenerate idempotent factors if and only if there exist an element
  \ $a\in G$ \ and a triangular array \ $\{\mu_{n,k}:n\in\NN,k=1,\ldots,K_n\}$ \
  of probability measures on \ $G$ \ such that
  \renewcommand{\labelenumi}{{\rm(\roman{enumi})}}
 \begin{enumerate}
   \item for every compact subset \ $C$ \ of \ $\hG,$ \
         $$
           \lim_{n\to\infty}\max_{1\sleq k\sleq K_n}\sup_{\chi\in C}\vert\hmu_{n,k}(\chi)-1\vert=0,
         $$
  \item for all \ $\chi\in\hG,$ \
        $$
         \sup_{n\in\NN}\sum_{k=1}^{K_n}(1-\vert\hmu_{n,k}(\chi)\vert)<+\infty,
        $$
  \item $\delta_a*\AstkoneKn\mu_{n,k}\stackrel{\mathrm{w}}{\longrightarrow}\mu$ \ as \ $n\to\infty,$ \
        where \ $\stackrel{\mathrm{w}}{\longrightarrow}$ \ means weak convergence and
        \ $\AstkoneKn\mu_{n,k}$ \
        denotes the convolution of \ $\mu_{n,k},$ $k=1,\ldots,K_n.$
 \end{enumerate}
\end{Thm}

We note that condition (i) of Theorem \ref{THM_YASUDA} is equivalent to the
infinitesimality of the triangular array \ $\{\mu_{n,k}:n\in\NN,k=1,\ldots,K_n\},$ \
see, e.g., 5.1.12 in Heyer \cite{HEY}.

\section{Weakly infinitely divisible measures on the torus}
\label{widT}

Consider the set \ $\TT:=\{\ee^{ix}:-\pi\leq x<\pi\}$ \ of complex numbers
under multiplication. This is a compact Abelian \ $T_0$--topological group having a
countable basis of its topology, and it is called the 1--dimensional torus group.
For elementary facts about \ $\TT$ \ we refer to the
monographs Hewitt and Ross \cite{HR}, Heyer \cite{HEY} and
 Hofmann and Morris \cite{HOFMOR}.
 The character group of \ $\TT$ \ is \ $\hTT=\{\chi_\ell:\ell\in\ZZ\}$, \ where
 $$
  \chi_\ell(y):=y^\ell,\qquad y\in\TT,\quad\ell\in\ZZ.
 $$
Hence \ $\hTT\cong\ZZ$ \ (i.e., \ $\hTT$ \ and \ $\ZZ$ \ are
topologically isomorphic).
The compact subgroups of \ $\TT$ \ are
 $$
  H_r:=\{\ee^{2\pi ij/r}:j=0,1,\dots,r-1\},\qquad r\in\NN,
 $$
 and \ $\TT$ \ itself.

The set of all quadratic forms on \ $\hTT\cong\ZZ$ \ is
 \ $\qq_+\big(\hTT\big)=\{\psi_b:b\in\RR_+\}$, \ where
 $$
  \psi_b(\chi_\ell):=b\ell^2,\qquad\ell\in\ZZ,\quad b\in\RR_+.
 $$
Let us define the functions \ $\arg:\TT\to[-\pi,\pi[$ \ and \
$h:\RR\to\RR$ \ by
 \begin{align*}
  \arg(\ee^{ix})&:=x,\qquad-\pi\leq x<\pi,\\
  h(x)&:=\begin{cases}
               0 & \text{if \ $x<-\pi$ \ or \ $x\geq\pi$,}\\
          -x-\pi & \text{if \ $-\pi\leq x<-\pi/2$,}\\
               x & \text{if \ $-\pi/2\leq x<\pi/2$,}\\
          -x+\pi & \text{if \ $\pi/2\leq x<\pi$.}\\
         \end{cases}
 \end{align*}
A measure \ $\eta$ \ on \ $\TT$ \ with values in \ $[0,+\infty]$ \ is a
 L\'evy measure if and only if \ $\eta(\{e\})=0$ \ and
 \ $\int_\TT(\arg y)^2\,\dd\eta(y)<\infty$.
\ The function \ $g_\TT:\TT\times\hTT\to\RR$,
 $$
  g_\TT(y,\chi_\ell):=\ell h(\arg y),\qquad y\in\TT,\quad\ell\in\ZZ,
 $$
 is a local inner product for \ $\TT$. \

Note that \ $\mathcal{I}(\TT)=\mathcal{I}_{\ww}(\TT),$ \
 since \ $(\ee^{ix/n})^n=\ee^{ix},$ $x\in[-\pi,\pi),$ $n\in\NN.$ \

Our aim is to show that for a weakly infinitely divisible measure \
$\mu$ \ on \ $\TT$ \ there exist independent real random variables \ $U$ \ and \ $Z$
 \ such that \ $U$ \ is uniformly distributed on a suitable subset of \ $\RR$,
 \ $Z$ \ has an infinitely divisible distribution on \ $\RR$, \ and
 \ $\ee^{i(U+Z)}\distre\mu$. \ We note that \ $\RR$ \ is a locally compact
   Abelian \ $T_0$--topological group, its character group is
     \ $\hRR=\{\chi_y:y\in\RR\}$, \ where \ $\chi_y(x):=\ee^{iyx}$. \
 The function \ $g_\RR:\RR\times\hRR\to\RR$, \ defined by
 \ $g_\RR(x,\chi_y):=yh(x)$, \ is a local inner product for \ $\RR$.

\begin{Thm} \
If \ $(H,a,\psi_b,\eta)\in\cP(\TT)$ \ then
 $$
  \ee^{i(U+\arg a+X+Y)}
  \distre\omega_H*\delta_a*\gamma_{\psi_b}*\pi_{\eta,\,g_\TT},
 $$
 where \ $U$, $X$ \ and \ $Y$ \ are independent real random variables such that
 \ $U$ \ is uniformly distributed on \ $[0,2\pi]$ \ if \ $H=\TT$, \ $U$ \ is
 uniformly distributed on \ $\{2\pi j/r:j=0,1,\dots,r-1\}$ \ if \ $H=H_r$,
 \ $X$ \ has a normal distribution on \ $\RR$ \ with zero mean and variance
 \ $b$, \ and the distribution of \ $Y$ \ is the generalized Poisson measure
 \ $\pi_{\arg\!\circ\eta,\,g_\RR}$ \ on \ $\RR$, \ where the measure
 \ $\arg\!\circ\eta$ \ on \ $\RR$ \ is defined by
 \ $(\arg\!\circ\eta)(B):=\eta\big(\{x\in\TT:\arg(x)\in B\}\big)$ \ for all
 Borel subsets \ $B$ \ of \ $\RR$.
\end{Thm}

\noindent{\bf Proof.} Let \ $U$ \ be a real random variable which is
uniformly distributed on
 \ $[0,2\pi]$.
\ Then for all \ $\chi_\ell\in\hTT$, \ $\ell\in\ZZ$, \ $\ell\not=0$,
 $$
  \EE\,\chi_\ell(\ee^{iU})
  =\EE\,\ee^{i\ell U}
  =\frac{1}{2\pi}\int_0^{2\pi}\ee^{i\ell x}\,\dd x
  =0.
 $$
Hence \ $\EE\,\chi_\ell(\ee^{iU})=\homega_\TT(\chi_\ell)$ \ for all
 \ $\chi_\ell\in\hTT$, \ $\ell\in\ZZ$, \ and we obtain
 \ $\ee^{iU}\distre\omega_\TT$.

Now let \ $U$ \ be a real random variable which is uniformly
distributed on
 \ $\{2\pi j/r:j=0,1,\dots,r-1\}$ \ with some \ $r\in\NN$.
\ Then for all \ $\chi_\ell\in\hTT$, \ $\ell\in\ZZ$,
 $$
  \EE\,\chi_\ell(\ee^{iU})
  =\EE\,\ee^{i\ell U}
  =\frac{1}{r}\sum_{j=0}^{r-1}\ee^{2\pi i\ell j/r}
  =\begin{cases}
    1 & \text{if \ $r|\ell$,}\\
    0 & \text{otherwise.}
   \end{cases}
 $$
Hence \ $\EE\,\chi_\ell(\ee^{iU})=\homega_{H_r}(\chi_\ell)$ \ for
all
 \ $\chi_\ell\in\hTT$, \ $\ell\in\ZZ$, \ and we obtain
 \ $\ee^{iU}\distre\omega_{H_r}$.

For \ $a\in\TT$, \ we have \ $a=\ee^{i\arg a}$, \ hence
 \ $\ee^{i\arg a}\distre\delta_a$.

For \ $b\in\RR_+$, \ the Fourier transform of the Gauss measure
 \ $\gamma_{\psi_b}$ \ has the form
 $$
  \hgamma_{\psi_b}(\chi_\ell)=\ee^{-b\ell^2/2},\qquad
  \chi_\ell\in\hTT,\quad\ell\in\ZZ.
 $$
For all \ $\chi_\ell\in\hTT$, \ $\ell\in\ZZ$,
 $$
  \EE\,\chi_\ell(\ee^{iX})
  =\EE\,\ee^{i\ell X}
  =\ee^{-b\ell^2/2}.
 $$
Hence \ $\EE\,\chi_\ell(\ee^{iX})=\gamma_{\psi_b}(\chi_\ell)$ \ for
all
 \ $\chi_\ell\in\hTT$, \ $\ell\in\ZZ$, \ and we obtain
 \ $\ee^{iX}\distre\gamma_{\psi_b}$.

For a L\'evy measure\ $\eta\in\LL(\TT)$, \ the Fourier transform of
the
 generalized Poisson measure \ $\pi_{\eta,\,g_\TT}$ \ has the form
 $$
  \hpi_{\eta,\,g_\TT}(\chi_\ell)
  =\exp\left\{\int_\TT\big(y^\ell-1-i\ell h(\arg y)\big)\dd\eta(y)\right\},
  \qquad\chi_\ell\in\hTT,\quad\ell\in\ZZ.
 $$
A measure \ $\widetilde{\eta}$ \ on \ $\RR$ \
with values in \ $[0,+\infty]$ \
is a L\'evy measure if and only if \ $\widetilde{\eta}(\{0\})=0$ \ and
 \ $\int_\RR\min\{1,x^2\}\,\dd\widetilde{\eta}(x)<\infty$.
\ Consequently, \ $\arg\!\circ\eta$ \ is a L\'evy measure on \
$\RR$, \ and for
 all \ $\chi_\ell\in\hTT$, \ $\ell\in\ZZ$,
 \begin{align*}
  \EE\,\chi_\ell(\ee^{iY})
  =\EE\,\ee^{i\ell Y}
  &=\exp\left\{\int_\RR
                \big(\ee^{i\ell x}-1-i\ell h(x)\big)
                 \,\dd(\arg\!\circ\eta)(x)\right\}\\[2mm]
  &=\exp\left\{\int_\TT\big(y^\ell-1-i\ell h(\arg y)\big)\,\dd\eta(y)\right\}.
 \end{align*}
Hence \ $\EE\,\chi_\ell(\ee^{iY})=\hpi_{\eta,\,g_\TT}(\chi_\ell)$ \
for all
 \ $\chi_\ell\in\hTT$, \ $\ell\in\ZZ$, \ and we obtain
 \ $\ee^{iY}\distre\pi_{\eta,\,g_\TT}$.

Finally, independence of \ $U$, $X$ \ and \ $Y$ \ implies
 \begin{align*}
  &\EE\,\chi(\ee^{i(U+\arg a+X+Y)})
   =\EE\,\chi(\ee^{iU})\cdot\chi(\ee^{i\arg a})\cdot\EE\,\chi(\ee^{iX})
    \cdot\EE\,\chi(\ee^{iY})\\
  &=\homega_H(\chi)\,\hdelta_a(\chi)\,\hgamma_{\psi_b}(\chi)\,
    \hpi_{\eta,\,g_\TT}(\chi)
   =(\omega_H*\delta_a*\gamma_{\psi_b}*\pi_{\eta,\,g_\TT})\:\widehat{}\:(\chi)
 \end{align*}
 for all \ $\chi\in\hTT$, \ hence we obtain the statement.
\proofend

\section{Weakly infinitely divisible measures on the group of \ $p$--adic
          integers}
\label{widD}

Let \ $p$ \ be a prime. The group of \ $p$--adic integers is
 $$
  \Delta_p
  :=\big\{(x_0,x_1,\dots)
          :\text{$x_j\in\{0,1,\dots,p-1\}$ \ for all \ $j\in\ZZ_+$}\big\},
 $$
 where the sum \ $z:=x+y\in\Delta_p$ \ for \ $x,y\in\Delta_p$ \ is uniquely
 determined by the relationships
 $$
  \sum_{j=0}^dz_jp^j\equiv\sum_{j=0}^d(x_j+y_j)p^j\quad\mod p^{d+1}
  \qquad\text{for all \ $d\in\ZZ_+$.}
 $$
For each \ $r\in\ZZ_+$, \ let
 $$
  \Lambda_r:=\{x\in\Delta_p:\text{$x_j=0$ \ for all \ $j\leq r-1$}\}.
 $$
The family of sets \ $\{x+\Lambda_r:x\in\Delta_p,\,r\in\ZZ_+\}$ \ is
an open subbasis for a topology on \ $\Delta_p$ \ under which \ $\Delta_p$ \ is a
 compact, totally disconnected Abelian \ $T_0$--topological group having a
 countable basis of its topology. For elementary facts about \ $\Delta_p$ \
 we refer to the monographs Hewitt and Ross \cite{HR}, Heyer \cite{HEY} and
 Hofmann and Morris \cite{HOFMOR}. The character group of \ $\Delta_p$ \ is
 \ $\hDelta_p=\{\chi_{d,\ell}:d\in\ZZ_+,\,\ell=0,1,\dots,p^{d+1}-1\}$, \ where
 $$
  \chi_{d,\ell}(x):=\ee^{2\pi i\ell(x_0+px_1+\cdots+p^dx_d)/p^{d+1}},\quad
  x\in\Delta_p,\quad d\in\ZZ_+,\;\ell=0,1,\dots,p^{d+1}-1.
 $$
The compact subgroups of \ $\Delta_p$ \ are \ $\Lambda_r$, \
$r\in\ZZ_+$
 \ (see Hewitt and Ross \cite[Example 10.16 (a)]{HR}).

 A measure \ $\eta$ \ on \ $\Delta_p$ \ with values in \ $[0,+\infty]$ \ is a
 L\'evy measure if and only if \ $\eta(\{e\})=0$ \ and
 \ $\eta(\Delta_p\setminus\Lambda_r)<\infty$ \ for all \ $r\in\ZZ_+$.

Since the group \ $\Delta_p$ \ is totally disconnected, the only
quadratic
 form on \ $\hDelta_p$ \ is \ $\psi=0$, \ and thus there is no
 nontrivial Gauss measure on \ $\Delta_p.$ \ Moreover, the function
 \ $g_{\Delta_p}:\Delta_p\times\hDelta_p\to\RR$, \ $g_{\Delta_p}=0$ \ is a
 local inner product for \ $\Delta_p$.

Now we prove that \ $\mathcal{I}(\Delta_p)\ne\mathcal{I}_{\ww}(\Delta_p)$ \
by showing that there exists an element \ $x\in\Delta_p$ \ such that
 \ $\delta_x\not\in\mathcal{I}(\Delta_p).$ \ Indeed, the Dirac measure at the element
 \ $(1,0,\dots)\in\Delta_p$ \ is not infinitely divisible, since for each element
  \ $y\in\Delta_p,$ \ the sum \ $py$ \ has the form \ $(0,z_1,z_2,\ldots)\in\Delta_p$ \
  with some \ $z_i\in\{0,1,\ldots,p-1\},$ $i\in\NN.$

Our aim is to show that for a weakly infinitely divisible measure \
$\mu$ \ on \ $\Delta_p$ \ there exist integer--valued random variables \ $U_0,U_1,\dots$
 \ and \ $Z_0,Z_1,\dots$ \ such that \ $U_0,U_1,\dots$ \ are independent of
 each other and of the sequence \ $Z_0,Z_1,\dots$, \ moreover \ $U_0,U_1,\dots$
 \ are uniformly distributed on a suitable subset of \ $\ZZ$,
 \ $(Z_0,\dots,Z_n)$ \ has a weakly infinitely divisible distribution on
 \ $\ZZ^{n+1}$ \ for all \ $n\in\ZZ_+$, \ and
 \ $\varphi(U_0+Z_0,\,U_1+Z_1,\dots)\distre\mu$, \ where the mapping
 \ $\varphi:\ZZ^\infty\to\Delta_p$, \ uniquely defined by the relationships
 \begin{equation}\label{cong}
  \sum_{j=0}^dy_jp^j\equiv\sum_{j=0}^d\varphi(y)_jp^j\quad\mod p^{d+1}
  \qquad\text{for all \ $d\in\ZZ_+$,}
 \end{equation}
 is a continuous homomorphism from the Abelian topological group \ $\ZZ^\infty$
 \ (furnished with the product topology) onto \ $\Delta_p$.
\ (Note that \ $\ZZ^\infty$ \ is not locally compact.) Continuity
of \ $\varphi$ \ follows from the definition of the product topology and the fact that
 $$
  \varphi^{-1}(x+\Lambda_r)
  =\{y\in\ZZ^\infty:(y_0,y_1,\dots,y_{r-1})\in F_{x,r}\}
 $$
 for all \ $x\in\Delta_p$, \ $r\in\ZZ_+$, \ where \ $F_{x,r}$ \ is a suitable (open)
 subset of \ $\ZZ^r$.

\begin{Thm}\label{repD} \
If \ $(\Lambda_r,a,0,\eta)\in\cP(\Delta_p)$ \ then
 $$
  \varphi(U_0+a_0+Y_0,\,U_1+a_1+Y_1,\dots)
  \distre\omega_{\Lambda_r}*\delta_a*\pi_{\eta,\,g_{\Delta_p}},
 $$
 where \ $U_0$, $U_1$, \dots \ and \ $Y_0$, $Y_1$, \dots \ are integer--valued
 random variables such that \ $U_0$, $U_1$, \dots \ are independent of each
 other and of the sequence \ $Y_0$, $Y_1$, \dots, \ moreover
 \ $U_0=\dots=U_{r-1}=0$ \ and \ $U_r$, $U_{r+1}$, \dots \ are uniformly
 distributed on \ $\{0,1,\dots,p-1\}$, \ and the distribution of
 \ $(Y_0,\dots,Y_n)$ \ is the compound Poisson measure \ $\ee(\eta_{n+1})$
 \ for all \ $n\in\ZZ_+$, \ where the measure \ $\eta_{n+1}$ \ on \ $\ZZ^{n+1}$
 \ is defined by \ $\eta_{n+1}(\{0\}):=0$ \ and
 \ $\eta_{n+1}(\ell):=\eta(\{x\in\Delta_p:(x_0,x_1,\dots,x_n)=\ell\})$ \ for
 all \ $\ell\in\ZZ^{n+1}\setminus\{0\}$.
\end{Thm}

\noindent{\bf Proof.} Since \ $U_0$, $U_1$, \dots \ and \ $Y_0$,
$Y_1$, \dots \ are integer--valued random variables and the mapping
  \ $\varphi:\ZZ^\infty\to\Delta_p$ \ is continuous,
we obtain that \ $\varphi(U_0+a_0+Y_0,\,U_1+a_1+Y_1,\dots)$ \ is a
 random element with values in \ $\Delta_p$.

First we show \ $\varphi(U)\distre\omega_{\Lambda_r}$, \ where
 \ $U:=(U_0,U_1,\dots)$.
\ By \eqref{cong} we obtain
 \begin{align*}
  &\EE\,\chi_{d,\ell}(\varphi(U))
   =\EE\,\ee^{2\pi i\ell
              (\varphi(U)_0+p\varphi(U)_1+\cdots+p^d\varphi(U)_d)/p^{d+1}}
   =\EE\,\ee^{2\pi i\ell(U_0+pU_1+\cdots+p^dU_d)/p^{d+1}}\\[2mm]
  &=\begin{cases}
     \DS\frac{1}{p^{d-r+1}}
      \sum_{j_r=0}^{p-1}\ldots\sum_{j_d=0}^{p-1}
       \ee^{2\pi i\ell(p^rj_r+\cdots+p^dj_d)/p^{d+1}}
     =0 & \text{if \ $d\geq r$ \ and \ $p^{d+1-r}\not|\ell$,}\\
     1 & \text{otherwise}\\
    \end{cases}
 \end{align*}
 for all \ $d\in\ZZ_+$ \ and \ $\ell=0,1,\dots,p^{d+1}-1$.
\ Hence \
$\EE\,\chi_{d,\ell}(\varphi(U))=\homega_{\Lambda_r}(\chi_{d,\ell})$
 \ for all \ $d\in\ZZ_+$ \ and \ $\ell=0,1,\dots,p^{d+1}-1$, \ and we obtain
 \ $\varphi(U)\distre\omega_{\Lambda_r}$.

For \ $a\in\Delta_p$, \ we have \ $a=\varphi(a_0,a_1,\dots)$, \
hence \ $\varphi(a_0,a_1,\dots)\distre\delta_a$.

For a L\'evy measure\ $\eta\in\LL(\Delta_p)$, \ the Fourier
transform of the generalized Poisson measure
  \ $\pi_{\eta,\,g_{\Delta_p}}$ \ has the form
 $$
  \hpi_{\eta,\,g_{\Delta_p}}(\chi_{d,\ell})
  =\exp\left\{\int_{\Delta_p}
    \big(\ee^{2\pi i\ell(x_0+px_1+\cdots+p^dx_d)/p^{d+1}}-1\big)
     \dd\eta(x)\right\}
 $$
 for all \ $d\in\ZZ_+$ \ and \ $\ell=0,1,\dots,p^{d+1}-1$.
\ We have
 \ $\eta_{n+1}(\ZZ^{n+1})=\eta(\Delta_p\setminus\Lambda_{n+1})<\infty$, \ hence
 \ $\eta_{n+1}$ \ is a bounded measure on \ $\ZZ^{n+1}$, \ and the compound
 Poisson measure \ $\ee(\eta_{n+1})$ \ on \ $\ZZ^{n+1}$ \ is defined.
The character group of \ $\ZZ^{n+1}$ \ is
 \ $(\ZZ^{n+1})\;\widehat{}\:
    =\{\chi_{z_0,z_1,\dots,z_n}:z_0,z_1,\dots,z_n\in\TT\}$,
 \ where
 \ $\chi_{z_0,z_1,\dots,z_n}(\ell_0,\ell_1,\dots,\ell_n)
    :=z_0^{\ell_0}z_1^{\ell_1}\cdots z_n^{\ell_n}$
 \ for all \ $(\ell_0,\ell_1,\dots,\ell_n)\in\ZZ^{n+1}$.
\ The family of measures \ $\{\ee(\eta_{n+1}):n\in\ZZ_+\}$ \ is
compatible,
 since \ $\ee(\eta_{n+2})(\{\ell\}\times\ZZ)=\ee(\eta_{n+1})(\{\ell\})$ \ for
 all \ $\ell\in\ZZ^{n+1}$ \ and \ $n\in\ZZ_+$.
\ Indeed, this is a consequence of
 $$
  (\ee(\eta_{n+2}))\:\widehat{}\:(\chi_{z_0,z_1,\dots,z_n,1})
  =(\ee(\eta_{n+1}))\:\widehat{}\:(\chi_{z_0,z_1,\dots,z_n})
 $$
 for all \ $z_0,z_1,\dots,z_n\in\TT$, \ which follows from
 \begin{align*}
  \int_{\ZZ^{n+2}}
   (z_0^{\ell_0}z_1^{\ell_1}\cdots z_n^{\ell_n}-1)\,
   &\dd\eta_{n+2}(\ell_0,\ell_1,\dots,\ell_n,\ell_{n+1})\\
   &=\int_{\ZZ^{n+1}}
     (z_0^{\ell_0}z_1^{\ell_1}\cdots z_n^{\ell_n}-1)\,
     \dd\eta_{n+1}(\ell_0,\ell_1,\dots,\ell_n)
 \end{align*}
 for all \ $z_0,z_1,\dots,z_n\in\TT$, \ where both sides are equal to
 \ $\int_{\Delta_p}(z_0^{x_0}z_1^{x_1}\cdots z_n^{x_n}-1)\,\dd\eta(x).$
 \ This integral is finite, since
 \begin{align*}
  \int_{\Delta_p}|z_0^{x_0}z_1^{x_1}\cdots z_n^{x_n}-1|\,\dd\eta(x)
  &=\int_{\Delta_p\setminus\Lambda_{n+1}}
    |z_0^{x_0}z_1^{x_1}\cdots z_n^{x_n}-1|\,\dd\eta(x)\\
  &\leq2\eta(\Delta_p\setminus\Lambda_{n+1})
   <\infty.
 \end{align*}
By Kolmogorov's Consistency Theorem, there exists a sequence \
$Y_0$, $Y_1$,
 \dots \ of integer--valued random variables such that the distribution of
 \ $(Y_0,\dots,Y_n)$ \ is the compound Poisson measure \ $\ee(\eta_{n+1})$
 \ for all \ $n\in\ZZ_+$.
\ For all \ $d\in\ZZ_+$ \ and \ $\ell=0,1,\dots,p^{d+1}-1$ \ we have
 \begin{align*}
  \EE\,\chi_{d,\ell}&(\varphi(Y_0,Y_1,\dots))
   =\EE\,\ee^{2\pi i\ell(Y_0+pY_1+\cdots+p^dY_d)/p^{d+1}}\\[2mm]
  &=\exp\left\{\int_{\ZZ^{d+1}}
     \big(\ee^{2\pi i\ell(\ell_0+p\ell_1+\cdots+p^d\ell_d)/p^{d+1}}-1\big)\,
      \dd\eta_{d+1}(\ell_0,\ell_1,\dots,\ell_d)\right\}\\[2mm]
  &=\exp\left\{\int_{\Delta_p}
     \big(\ee^{2\pi i\ell(x_0+px_1+\cdots+p^dx_d)/p^{d+1}}-1\big)\,
      \dd\eta(x)\right\}.
 \end{align*}
Hence
 \ $\EE\,\chi_{d,\ell}(\varphi(Y_0,Y_1,\dots))
    =\hpi_{\eta,\,g_{\Delta_p}}(\chi_{d,\ell})$
 \ for all \ $d\in\ZZ_+$ \ and \ $\ell=0,1,\dots,p^{d+1}-1$, \ and we obtain
 \ $\varphi(Y_0,Y_1,\dots)\distre\pi_{\eta,\,g_{\Delta_p}}$.

Since the sequences \ $U_0,U_1,\dots$ \ and \ $Y_0,Y_1,\dots$ \ are
 independent and the mapping \ $\varphi:\ZZ^\infty\to\Delta_p$ \ is a
 homomorphism, we have
 \begin{align*}
  &\EE\,\chi(\varphi(U_0+a_0+Y_0,\,U_1+a_1+Y_1,\dots))\\
  &=\EE\,\chi(\varphi(U_0,U_1,\dots))\cdot
    \chi(\varphi(a_0,a_1,\dots))
    \cdot\EE\,\chi(\varphi(Y_0,Y_1,\dots))\\
  &=\homega_{\Lambda_r}(\chi)\,\hdelta_a(\chi)\,
    \hpi_{\eta,\,g_{\Delta_p}}(\chi)
   =(\omega_{\Lambda_r}*\delta_a*\pi_{\eta,\,g_{\Delta_p}})\:\widehat{}\:(\chi)
 \end{align*}
 for all \ $\chi\in\hDelta_p$, \ and we obtain the statement.
\proofend

\section{Weakly infinitely divisible measures on the \ $p$--adic solenoid}
\label{widS}

Let \ $p$ \ be a prime. The \ $p$--adic solenoid is a subgroup of \
$\TT^\infty$, \ namely,
 $$
  S_p=\big\{(y_0,y_1,\dots)\in\TT^\infty:
            \text{$y_j=y_{j+1}^p$ \ for all \ $j\in\ZZ_+$}\big\}.
 $$
This is a compact Abelian \ $T_0$--topological group having a
countable basis of its topology. For elementary facts about \ $S_p$ \
 we refer to the monographs Hewitt and Ross \cite{HR}, Heyer \cite{HEY} and
 Hofmann and Morris \cite{HOFMOR}.
The character group \ $S_p$ \ is
 \ $\hS_p=\{\chi_{d,\ell}:d\in\ZZ_+,\,\ell\in\ZZ\}$, \ where
 $$
  \chi_{d,\ell}(y):=y_d^\ell,\qquad
  y\in S_p,\quad d\in\ZZ_+,\quad\ell\in\ZZ.
 $$
The set of all quadratic forms on \ $\hS_p$ \ is
 \ $\qq_+\big(\hS_p\big)=\{\psi_b:b\in\RR_+\}$, \ where
 $$
  \psi_b(\chi_{d,\ell}):=\frac{b\ell^2}{p^{2d}},\qquad
  d\in\ZZ_+,\quad\ell\in\ZZ,\quad b\in\RR_+.
 $$
A measure \ $\eta$ \ on \ $S_p$ \ with values in \ $[0,+\infty]$ \ is a
 L\'evy measure if and only if \ $\eta(\{e\})=0$ \ and
 \ $\int_{S_p}(\arg y_0)^2\,\dd\eta(y)<\infty$.
\ The function \ $g_{S_p}:S_p\times\hS_p\to\RR$,
 $$
  g_{S_p}(y,\chi_{d,\ell}):=\frac{\ell h(\arg y_0)}{p^d},\qquad
  y\in S_p,\quad d\in\ZZ_+,\quad\ell\in\ZZ,
 $$
 is a local inner product for \ $S_p$.

Using a result of Carnal \cite{CAR}, Becker-Kern \cite{BEC} showed that
all the Dirac measures on \ $S_p$ \ are infinitely divisible, which implies 
 \ $\mathcal{I}(S_p)=\mathcal{I}_{\ww}(S_p).$ \

Our aim is to show that for a weakly infinitely divisible measure \
$\mu$ \ on
 \ $S_p$ \ without an idempotent factor there exist real random variables
 \ $Z_0,Z_1,\dots$ \ such that \ $(Z_0,\dots,Z_n)$ \ has a weakly infinitely
 divisible distribution on \ $\RR\times\ZZ^n$ \ for all \ $n\in\ZZ_+$, \ and
 \ $\varphi(Z_0,\,Z_1,\dots)\distre\mu$, \ where the mapping
 \ $\varphi:\RR\times\ZZ^\infty\to S_p$, \ defined by
 \begin{align*}
  \varphi&(y_0,y_1,y_2,\dots)\\
    &:=\big(\ee^{iy_0},\,\ee^{i(y_0+2\pi y_1)/p},\,
         \ee^{i(y_0+2\pi y_1+2\pi y_2p)/p^2},\,
         \ee^{i(y_0+2\pi y_1+2\pi y_2p+2\pi y_3p^2)/p^3},\dots\big)
 \end{align*}
 for \ $(y_0,y_1,y_2,\dots)\in\RR\times\ZZ^\infty$, \ is a continuous
 homomorphism from the Abelian topological group \ $\RR\times\ZZ^\infty$
 \ (furnished with the product topology) onto \ $S_p$. \
 Continuity of \ $\varphi$ \ follows from the fact that 
  \ $S_p$ \ as a subspace of \ $\TT^\infty$ \ is furnished with the relative topology. 
 Note that \ $\RR\times\ZZ^\infty$ \ is not locally compact, but
 \ $\RR\times\ZZ^n$ \ is a locally compact Abelian \ $T_0$--topological group
 having a countable basis of its topology for all \ $n\in\ZZ_+$.
\ The character group of \ $\RR\times\ZZ^n$ \ is
 \ $(\RR\times\ZZ^n)\:\widehat{}\:=\{\chi_{y,z}:y\in\RR,\,z\in\TT^n\}$, \ where
 \ $\chi_{y,z}(x,\ell):=\ee^{iyx}z_1^{\ell_1}\cdots z_n^{\ell_n}$ \ for all
 \ $x,y\in\RR$, \ $z=(z_1,\dots,z_n)\in\TT^n$ \ and
 \ $\ell=(\ell_1,\dots,\ell_n)\in\ZZ^n$.
\ The function \
$g_{\RR\times\ZZ^n}\big((x,\ell),\chi_{y,z}\big):=yh(x)$ \ is
 a local inner product for \ $\RR\times\ZZ^n$.

We also find independent real random variables $U_0,U_1,\dots$ such
that $U_0,U_1,\dots$ \ are uniformly distributed on suitable subsets
of \ $\RR$ \
 and \ $\varphi(U_0,\,U_1,\dots)\distre\omega_{S_p}$.

\begin{Thm}\label{TwidS} \
If \ $(\{e\},a,\psi_b,\eta)\in\cP(S_p)$ \ then
 \begin{align*}
  &\varphi(\tau(a)_0+X_0+Y_0,\,\tau(a)_1+Y_1,\,\tau(a)_2+Y_2,\dots)\\
  &=\Big(a_0\ee^{i(X_0+Y_0)},\,a_1\ee^{i(X_0+Y_0+2\pi Y_1)/p},\,
         a_2\ee^{i(X_0+Y_0+2\pi Y_1+2\pi Y_2p)/p^2},\dots\Big)\\
  &\distre\delta_a*\gamma_{\psi_b}*\pi_{\eta,\,g_{S_p}},
 \end{align*}
 where the mapping \ $\tau:S_p\to\RR\times\ZZ^\infty$ \ is defined by
 $$
   \tau(x)
         :=\left(\arg x_0,\,\frac{p\arg x_1-\arg x_0}{2\pi},\,
                 \frac{p\arg x_2-\arg x_1}{2\pi},\dots\right)
 $$
 for \ $x=(x_0,x_1,\dots)\in S_p.$ \ Here \ $X_0$, $Y_0$ \ are real random
 variables and \ $Y_1,Y_2,\dots$ \ are
 integer--valued random variables such that \ $X_0$ \ is independent of the
 sequence \ $Y_0,Y_1,\dots$, \ the variable \ $X_0$ \ has a normal distribution
 with zero mean and variance \ $b$, \ and the distribution of
 \ $(Y_0,\dots,Y_n)$ \ is the generalized Poisson measure
 \ $\pi_{\eta_{n+1},\,g_{\RR\times\ZZ^n}}$ \ for all \ $n\in\ZZ_+$,
 \ where the measure \ $\eta_{n+1}$ \ on \ $\RR\times\ZZ^n$ \ is defined by
 \ $\eta_{n+1}(\{0\}):=0$ \ and
 $$
  \eta_{n+1}(B\times\{\ell\})
  :=\eta\big(\big\{x\in S_p
                   :\tau(x)_0\in B,\,
                    (\tau(x)_1,\dots,\tau(x)_n)=\ell\big\}\big)
 $$
 for all Borel subsets \ $B$ \ of \ $\RR$ \ and for all \ $\ell\in\ZZ^n$ \ with
 \ $0\not\in B\times\{\ell\}$.

Moreover,
 $$
  \varphi(U_0,\,U_1,\dots)
  \distre\omega_{S_p},
 $$
 where \ $U_0,U_1,\dots$ \ are independent real random variables such that
 \ $U_0$ \ is uniformly distributed on \ $[0,2\pi]$ \ and
 \ $U_1,U_2,\dots$ \ are uniformly distributed on \ $\{0,1,\dots,p-1\}$.
\end{Thm}

\noindent{\bf Proof.} Since \ $X_0$, $Y_0$ \ and \ $U_0,U_1,\dots$ \
are real random variables and
 \ $Y_1,Y_2,\dots$ \ are integer--valued random variables and the mapping
 \ $\varphi:\RR\times\ZZ^\infty\to S_p$ \ is continuous, we obtain that
 \ $\varphi(\tau(a)_0+X_0+Y_0,\,\tau(a)_1+Y_1,\,\tau(a)_2+Y_2,\dots)$ \ and
 \ $\varphi(U_0,\,U_1,\dots)$ \ are random elements with values in \ $S_p$.

For \ $a\in S_p$, \ we have \ $a=\varphi(\tau(a))$, \ hence
 \ $\varphi(\tau(a))\distre\delta_a$.

For \ $b\in\RR_+$, \ the Fourier transform of the Gauss measure
 \ $\gamma_{\psi_b}$ \ has the form
 $$
  \hgamma_{\psi_b}(\chi_{d,\ell})
  =\exp\left\{-\frac{b\ell^2}{2p^{2d}}\right\},\qquad
  d\in\ZZ_+,\quad\ell\in\ZZ.
 $$
For all \ $d\in\ZZ_+$ \ and \ $\ell\in\ZZ$,
 $$
  \EE\,\chi_{d,\ell}(\varphi(X_0,0,0,\dots))
  =\EE\,\ee^{i\ell X_0/p^d}
  =\exp\left\{-\frac{b\ell^2}{2p^{2d}}\right\}.
 $$
Hence
 \ $\EE\,\chi_{d,\ell}(\varphi(X_0,0,0,\dots))=\hgamma_{\psi_b}(\chi_{d,\ell})$
 \ for all \ $d\in\ZZ_+$ \ and \ $\ell\in\ZZ$, \ and we obtain
 that \ $\varphi(X_0,0,0,\dots)\distre\gamma_{\psi_b}$.

For a L\'evy measure \ $\eta\in\LL(S_p)$, \ the Fourier transform of
the generalized Poisson measure \ $\pi_{\eta,\,g_{S_p}}$ \ has the
form
 $$
  \hpi_{\eta,\,g_{S_p}}(\chi_{d,\ell})
  =\exp\left\{\int_{S_p}
               \big(y_d^\ell-1-i\ell h(\arg y_0)/p^d\big)\,\dd\eta(y)\right\}
 $$
 for all \ $d\in\ZZ_+$ \ and \ $\ell\in\ZZ$. \
 A measure \ $\widetilde{\eta}$ \ on \ $\RR\times\ZZ^n$ \
 with values in \ $[0,+\infty]$ \ is a L\'evy measure
 if and only if \ $\widetilde{\eta}(\{0\})=0$,
 \ $\widetilde{\eta}(\{(x,\ell)\in\RR\times\ZZ^n
                       :\text{$|x|\geq\vare$ or $\ell\not=0$}\})<\infty$
 \ for all \ $\vare>0$, \ and
 \ $\int_{\RR\times\ZZ^n}h(x)^2\,\dd\widetilde{\eta}(x,\ell)<\infty$.
\ We have
 \begin{align*}
  &\eta_{n+1}(\{(x,\ell)\in\RR\times\ZZ^n
                :\text{$|x|\geq\vare$ or $\ell\not=0$}\})\\
  &=\eta(\{y\in S_p
           :\text{$|\arg y_0|\geq\vare$ or
                  $(\tau(y)_1,\dots,\tau(y)_n)\not=0$}\})
   =\eta(S_p\setminus N_{\vare,n})
    <\infty
 \end{align*}
 for all \ $\vare\in(0,\pi)$, \ where
 $$
  N_{\vare,n}
  :=\{y\in S_p
      :|\arg y_0|<\vare,\,|\arg y_1|<\vare/p,\dots,|\arg y_n|<\vare/p^n\}.
 $$
Moreover,
 \ $\int_{\RR\times\ZZ^n}h(x)^2\,\dd\eta_{n+1}(x,\ell)
    =\int_{S_p}h(\arg y_0)^2\,\dd\eta(y)<\infty$,
 \ since \ $\eta$ \ is a L\'evy measure on \ $S_p$.
\ Consequently, \ $\eta_{n+1}$ \ is a L\'evy measure on \
$\RR\times\ZZ^n$. \ The family of measures
 \ $\{\pi_{\eta_{n+1},\,g_{\RR\times\ZZ^n}}:n\in\ZZ_+\}$ \ is compatible, since
 \ $\pi_{\eta_{n+2},\,g_{\RR\times\ZZ^{n+1}}}(\{x\}\times\ZZ)
    =\pi_{\eta_{n+1},\,g_{\RR\times\ZZ^n}}(\{x\})$
 \ for all \ $x\in\RR\times\ZZ^{n+1}$ \ and \ $n\in\ZZ_+$.
\ Indeed, this is a consequence of
 $$
  (\pi_{\eta_{n+2},\,g_{\RR\times\ZZ^{n+1}}})\:\widehat{}\:
   (\chi_{y,z_1,\dots,z_n,1})
  =(\pi_{\eta_{n+1},\,g_{\RR\times\ZZ^n}})\:\widehat{}\:
    (\chi_{y,z_1,\dots,z_n})
 $$
 for all \ $y\in\RR$, \ $z_1,\dots,z_n\in\TT$, \ which follows from
 \begin{align*}
  &\int_{\RR\times\ZZ^{n+1}}
    \big(\ee^{iyx}z_1^{\ell_1}\cdots z_n^{\ell_n}-1-iyh(x)\big)\,
     \dd\eta_{n+2}(x,\ell_1,\dots,\ell_n,\ell_{n+1})\\[2mm]
  &=\int_{\RR\times\ZZ^n}
    \big(\ee^{iyx}z_1^{\ell_1}\cdots z_n^{\ell_n}-1-iyh(x)\big)\,
     \dd\eta_{n+1}(x,\ell_1,\dots,\ell_n)
 \end{align*}
 for all \ $y\in\RR$, \ $z_1,\dots,z_n\in\TT$, \ where both sides are equal to
 \begin{align*}
  I:=\int_{S_p}
      \big(&\ee^{iy\arg x_0}z_1^{(p\arg x_1-\arg x_0)/(2\pi)}\cdots
                z_n^{(p\arg x_n-\arg x_{n-1})/(2\pi)}\\
      &-1-iyh(\arg x_0)\big)\,\dd\eta(x).
 \end{align*}
This integral is finite. Indeed, for all \ $x\in N_{\vare,n}$ \ and
\ $0<\vare<\pi/2$ \ we have
 \ $p\arg x_k=\arg x_{k-1}$ \ for each \ $k=1,\dots,n$, \ hence
 \begin{align*}
  |I|&\leq(2+\pi|y|)\,\eta(S_p\setminus N_{\vare,n})
          +\int_{N_{\vare,n}}|\ee^{iy\arg x_0}-1-iy\arg x_0|\,\dd\eta(x)\\[2mm]
     &\leq(2+\pi|y|)\,\eta(S_p\setminus N_{\vare,n})
          +\frac{1}{2}\int_{N_{\vare,n}}(\arg x_0)^2\,\dd\eta(x)
     <\infty,
 \end{align*}
 since \ $\eta$ \ is a L\'evy measure on \ $S_p$.
\ By Kolmogorov's Consistency Theorem, there exist a real random
variable \ $Y_0$ \ and a sequence \ $Y_1$, $Y_2$, \dots \ of integer--valued random
 variables such that the distribution of \ $(Y_0,\dots,Y_n)$ \ is the
 generalized Poisson measure \ $\pi_{\eta_{n+1},\,g_{\RR\times\ZZ^n}}$ \ for
 all \ $n\in\ZZ_+$.
\ For all \ $d\in\ZZ_+$ \ and \ $\ell\in\ZZ,$ \ we have
 \begin{align*}
  &\EE\,\chi_{d,\ell}(\varphi(Y_0,Y_1,\dots))
   =\EE\,\ee^{i\ell(Y_0+2\pi Y_1+\cdots+2\pi Y_dp^{d-1})/p^d}\\[2mm]
  &=\exp\left\{\int_{\RR\times\ZZ^d}
                \big(\ee^{i\ell(x+2\pi\ell_1+\cdots+2\pi\ell_dp^{d-1})/p^d}
                     -1-i\ell h(x)/p^d\big)\,
                \dd\eta_{d+1}(x,\ell_1,\dots,\ell_d)\right\} \\[2mm]
  &=\exp\left\{\int_{S_p}
                \big(y_d^\ell-1-i\ell h(\arg y_0)/p^d\big)\,\dd\eta(y)\right\}.
 \end{align*}
Hence
 \ $\EE\,\chi_{d,\ell}(\varphi(Y_0,Y_1,\dots))
    =\hpi_{\eta,\,g_{S_p}}(\chi_{d,\ell})$
 \ for all \ $d\in\ZZ_+$ \ and \ $\ell\in\ZZ$, \ and we obtain
 \ $\varphi(Y_0,Y_1,\dots)\distre\pi_{\eta,\,g_{S_p}}$.

Since the sequence \ $Y_0,Y_1,\dots$ \ and the random variable \
$X_0$ \ are
 independent and the mapping \ $\varphi:\RR\times\ZZ^\infty\to S_p$ \ is a
 homomorphism, we get
 \begin{align*}
  &\EE\,\chi(\varphi(\tau(a)_0+X_0+Y_0,\,\tau(a)_1+Y_1,\,
             \tau(a)_2+Y_2,\dots))\\
  &=\chi(\varphi(\tau(a)_0,\tau(a)_1,\dots))
    \cdot\EE\,\chi(\varphi(X_0,0,0,\dots))
    \cdot\EE\,\chi(\varphi(Y_0,Y_1,\dots))\\
  &=\hdelta_a(\chi)\,\hgamma_{\psi_b}(\chi)\,\hpi_{\eta,\,g_{S_p}}(\chi)
   =(\delta_a*\gamma_{\psi_b}*\pi_{\eta,\,g_{S_p}})\:\widehat{}\:(\chi)
 \end{align*}
 for all \ $\chi\in\hS_p$, \ and we obtain the first statement.

For all \ $d\in\ZZ_+$ \ and \ $\ell\in\ZZ\setminus\{0\}$,
 \begin{align*}
  &\EE\,\chi_{d,\ell}(\varphi(U_0,U_1,\dots))
   =\EE\,\ee^{i\ell(U_0+2\pi U_1+\cdots+2\pi U_dp^{d-1})/p^d}\\[2mm]
  &=\frac{1}{2\pi p^d}\int_0^{2\pi}\ee^{i\ell x/p^d}\,\dd x
    \sum_{j_0=0}^{p-1}\ldots\sum_{j_{d-1}=0}^{p-1}
     \ee^{2\pi i\ell(j_0+j_1p+\cdots+j_{d-1}p^{d-1})/p^d}
   =0.
 \end{align*}
Hence
 \ $\EE\,\chi_{d,\ell}(\varphi(U_0,U_1,\dots))=\homega_{S_p}(\chi_{d,\ell})$
 \ for all \ $d\in\ZZ_+$ \ and \ $\ell\in\ZZ$, \ and we obtain
 \ $\varphi(U_0,U_1,\dots)\distre\omega_{S_p}$.
\proofend

\vskip0.2cm

\noindent{\bf \large Acknowledgments.} The authors have been
supported by the Hungarian Scientific Research Fund under Grant No.\
OTKA--T048544/2005. The first author has been also supported by the
Hungarian Scientific Research Fund under Grant No.\
OTKA--F046061/2004.

\end{document}